\documentclass{ifacconf}
\usepackage[authoryear,round]{natbib}
\usepackage{graphics} 
\usepackage{epsfig} 
\usepackage{mathptmx} 
\usepackage{fancyhdr}
\usepackage{times} 
\usepackage{amsmath} 
\usepackage{amssymb}  
\usepackage[utf8]{inputenc}
\usepackage[T1]{fontenc}
\usepackage[section]{placeins}
\usepackage[export]{adjustbox}
\usepackage{float}
\usepackage{lastpage}
\usepackage{url}
\usepackage{color}
\usepackage{graphicx}     
\usepackage{natbib}

\begin{document}
\begin{frontmatter}

\title{\emph{A posteriori} probabilistic feasibility guarantees for Nash equilibria in uncertain multi-agent games\thanksref{footnoteinfo}}

\thanks[footnoteinfo]{Research was supported by the UK Engineering and Physical Sciences Research Council (EPSRC) under grant agreement EP/P03277X/1. }

\author{George Pantazis, } 
\author{Filiberto Fele, } 
\author{Kostas Margellos}

\address{Department of Engineering Science, University of Oxford, OX1 3PJ, UK, (e-mail: georgios.pantazis@lmh.ox.ac.uk,  filiberto.fele@eng.ox.ac.uk,  kostas.margellos@eng.ox.ac.uk)}

\begin{abstract}                
In this paper a distribution-free methodology is presented for providing robustness guarantees for Nash equilibria (NE) of multi-agent games. Leveraging recent \emph{a posteriori} developments of the so called scenario approach \citep{CampiGarattiRamponi2018}, we provide probabilistic guarantees for feasibility problems with polytopic constraints. This result is then used in the context of multi-agent games, allowing to provide robustness certificates for  constraint violation of any NE of a given game. Our guarantees can be used alongside any NE seeking algorithm that returns some equilibrium solution. Finally, by exploiting the structure of our problem, we circumvent the need of employing computationally prohibitive algorithms to find an irreducible support subsample, a concept at the core of the scenario approach. Our theoretical results are accompanied by simulation studies that investigate the robustness of the solutions of two different problems, namely, a 2-dimensional feasibility problem and an electric vehicle (EV) charging control problem.
\end{abstract}
\begin{keyword}
Scenario approach, Multi-agent games, Nash equilibria, Feasibility guarantees, Electric vehicles
\end{keyword}

\end{frontmatter}

\section{Introduction} \label{sec:secI}

Decentralized optimization and control of large scale systems is of significant interest for a variety of fields  ranging from  engineering and biology to economics and social sciences.  In many cases, systems of this kind can be modelled as a network comprising self-interested entities/agents  that interact/compete with each other in order to meet their own individual goals, thus giving rise to a noncooperative set-up. Aiming at the resolution of the inherent conflict among the agents, a vast amount of research work  has been focused on game-theoretic approaches. Especially in engineering, this framework has been extensively used for the analysis of communication  \citep{Pang2014}, \citep{Basar2002} and traffic  \citep{Smith1979} networks, smart grids \citep{Saad2012}, \citep{Callaway2013} and  electricity markets \citep{Chen2014}.  \par
  In real-world applications, however, the system is affected by uncertainty that can be due to several factors. 
This gives rise to games of incomplete information \citep{HarsanyiGames}. The concept of Nash equilibrium (NE) (see Definition \ref{def1}) typically assumes complete information. For incomplete information games the notion of NE appears insufficient to deal with the presence of uncertainty, as the resulting agents' strategies do not necessarily exhibit any robustness properties against uncertainty. \par 
  One of the first attempts to deal with uncertainty in multi-agent games was noticed in \citep{Harsanyi1962}. 
Motivated by these developments, two main research directions are encountered: 1) NE analysis based on particular models for the probability distribution of the uncertainty  \citep{HarsanyiGames},  \citep{Kouvaritakis2005}, \citep{Singh2016}  and/or the geometry of its support \citep{Aghassi2006}, \citep{Fukushima2005} and, 2) Distribution-free NE analysis, where no assumption on the probability distribution of the uncertainty is imposed. \par

 In this paper we focus on distribution-free NE seeking where algorithmic developments have been quite restricted. Motivated by the lack of distribution-free results, we leverage  the recent developments of the so called scenario approach \citep{CampiGarattiRamponi2018}  and represent the uncertainty by means of scenarios that could be either available as historical data, or extracted via some prediction model. As such, the two main theoretical pillars of this work are game theory and the scenario approach. We attempt to transfer concepts from the scenario approach, well-understood in the context of optimization, to multi-agent games aiming to accompany game equilibria with guarantees on the probability of constraint satisfaction.  \par
The scenario approach is  a well-established mathematical technique  \citep{CampiCalafiore2006}, \citep{Campi2008a}, \citep{Campi2008},  and  still a highly active research area (see \citep{Campi2018}, \citep{CampiGarattiRamponi2018} for some recent results), 
originally introduced to provide \emph{a priori} probabilistic guarantees for solutions of uncertain convex optimization programs. Very recently, the theory has been extended to non-convex decision making problems \citep{CampiGarattiRamponi2018}, where the probabilistic guarantees are obtained in an \emph{a posteriori} fashion. The main advantage of the scenario approach is its applicability under very general conditions, since it does not require the knowledge of the uncertainty set or the probability distribution, key assumptions in robust \citep{Bai1997} and stochastic optimization \citep{Birge1997}, respectively.\par According to the scenario approach, the original problem can be approximated by solving a computationally tractable approximate problem, the so called scenario program consisting of a finite number of constraints, each of them corresponding to a different realization of the uncertain parameter. Apart from its simplicity, the important feature of this program is its generalization properties, i.e., using a relatively small number of samples, guarantees for the probability that its solution satisfies a yet unseen constraint is obtained.   \par

Only a few data-driven works for distribution-free NE seeking have appeared in the literature with \citep{Paccagnan2019a} and  \citep{Fele2019a}, \citep{Fele2019b} being the most closely related to our work. Both papers attempt to bridge multi-agent games with the scenario approach, thus providing a distribution-free way to determine NE with quantifiable robustness properties. Specifically, in \citep{Paccagnan2019a}, the authors focus on variational inequalities, a typical solution concept for multi-agent games affected by uncertainty represented by means of scenarios. Due to the connection between variational inequalities and NE, their developments naturally find applications to multi-agent games. However, the theoretical analysis requires strong monotonicity of the operator associated with the variational inequalities, which in a gaming setting implies uniqueness of the NE. The latter is quite restrictive, as many games of practical interest have multiple (possibly infinite) NE.
In \citep{Fele2019a} and \citep{Fele2019b}, multi-agent games with uncertainty affecting agents' cost functions are considered, and once again uncertainty is tackled by means of the scenario approach. However, agents' constraint sets are assumed to be deterministic. Uncertain constraints could still be considered, under the restrictive assumption that they are decoupled. The approach proposed here is instead applicable to general uncertain constraints. 
The contributions of our paper with respect to the aforementioned works are the following: 
\begin{enumerate}
\item Leveraging the recent results of \citep{CampiGarattiRamponi2018}, we provide \emph{a posteriori} robustness certificates for the entire feasibility region of feasibility programs with polytopic constraints.
\item Focusing on multi-agent games, we provide distribution-free probabilistic guarantees for the entire set of NE in an \emph{a posteriori} fashion. We extend the results of \citep{Fele2019a} to account for uncertainty in agents' (possibly) coupling constraints, while we do not require uniqueness of the associated NE as in \citep{Paccagnan2019a}.  
\item The probabilistic results of recent works in an \emph{a posteriori} context \citep{Campi2018},  
 \citep{CampiGarattiRamponi2018} rely on certain algorithms for quantifying a so called irreducible support subsample (see Definition 2 in \citep{Campi2018}). Restricting our attention to multi-agent games, the use of such algorithms alongside NE seeking iterative algorithms, apart from being computationally prohibitive, leads to erratic behaviour due to numerical issues. In our case, the cardinality of the support subsample coincides with the number of facets of the polytopic constraint set, which is directly available. 
\end{enumerate}

The rest of the paper is organized as follows: Section 2 introduces the problem under study and offers a motivating application. Section 3 provides the theoretical analysis and proof of the main results, namely, providing \emph{a posteriori} guarantees for feasibility problems. Section 4 provides numerical examples  and revisits the electric vehicle application of Section 2. Finally, Section 5 concludes the paper and provides some directions for future work.

\section{Scenario-based multi-agent games with uncertain constraints}
\subsection{Gaming set-up}

Let $M$ be the total number of agents and $x^m=(x^m_t)_{t=1}^d$ the decision vector of dimension $d \in \mathbb{N}$ of agent $m \in \mathcal{M}=\{1,...,M \}$ taking values in the uncertain set $X^m_\delta \subseteq \mathbb{R}^d$ parameterized by the uncertain parameter $\delta$. The uncertainty parameter is defined on the (possibly unknown) probability space $(\Delta,\mathcal{F}, \mathbb{P})$, where $\Delta$ is the sample space, equipped with a $\sigma$-algebra   $\mathcal{F}$  and a probability measure  $\mathbb{P}$.  
 Similarly, we define $x^{-m}=(x^1,...,x^{m-1},x^{m+1},...,x^M)\in X_\delta^{-m}$, where $X^{-m}_\delta=\prod_{j\in \mathcal{M}, j  \neq m} X^j_\delta \subseteq \mathbb{R}^{d(M-1)}$, as the vector comprising  the decisions of all other agents except for that of agent $m$. Finally, let $\{\delta_i\}_{i=1}^N \in \Delta^N$ be a finite collection of independent and identically distributed (i.i.d.) scenarios/realisations of the uncertain vector $\delta$, where $\Delta^N$ is the cartesian product of multiple copies of the sample space $\Delta$. 
 In our set-up, agents are considered as self-interested entities, i.e., they are interested in minimizing their own deterministic cost function $J_m: X  \rightarrow \mathbb{R}$, where $ \mathbb{X}=\prod_{m=1}^M X^m \subseteq \mathbb{R}^{dM}$. \par
 Furthermore, we impose the following assumption: 
\begin{assumption} \label{affine1}
  Each agent's decision set is formed by the intersection of a deterministic decision space $X^m$ and uncertain affine constraints affected  by any realization $\delta \in \Delta$, i.e., $ \bigcap_{\delta \in \Delta} X^m_\delta =\{x^m \in X^m : g(x^m,x^{-m},\delta) \leq 0 \}, \ \forall \ m \in \mathcal{M}$\footnote{Formally, the intersection $\bigcap_{\delta \in \Delta} X^m_\delta$ of agent $m$ possibly depends on the strategies of all other agents, thus allowing the application of our results in generalised NE problems. For simplicity we drop this dependence in the subsequent analysis, as our main focus is the treatment of uncertainty.}, where $g: \mathbb{X} \times \Delta \rightarrow \mathbb{R}$ is an affine function with respect to its first two arguments.
\end{assumption}

Under Assumption \ref{affine1}, we consider a multi-agent game, whose constraints are affected by uncertainty. Each agent $m \in \mathcal{M}$ seeks to minimize her own cost function, given the strategies $x^{-m}$ of all other players, by solving the following program
\begin{align}
 \min_{x^m \in X^m}{J_m(x^m,x^{-m})} \ \text{subject to} \  x^m  \in \bigcap_{\delta \in \Delta} X^m_\delta.  \label{EVprin} 
\end{align}
For the convenience of the reader, we recall that $\delta$ is an uncertain parameter defined on the (possibly) unknown probability space $(\Delta,\mathcal{F},\mathbb{P})$.
For the aforementioned problem, we consider the solution concept of NE as presented in  Definition \ref{def1}.  
\begin{definition} \label{def1} \citep{Basar1999}
A vector $x_{NE}=\\
(x^m_{NE})_{m \in \mathcal{M}}$ is a NE of the associated game if and only if $J_m(x^m_{NE},x^{-m}_{NE}) \leq J_m(x^m,x^{-m}_{NE}) $ for any $x^m \in  \bigcap_{\delta \in \Delta} X^m_\delta$ and for any $m \in \mathcal{M}$.
\end{definition}
 Due to the presence of uncertainty and the (possibly) infinite cardinality of $\Delta$, problem (\ref{EVprin}) is very difficult to solve, without imposing any assumptions on the geometry of the sample set $\Delta$ or the underlying probability distribution $\mathbb{P}$. 
To circumvent those issues, we approximate problem (\ref{EVprin}) by drawing multiple i.i.d. samples $\{\delta_i\}_{i=1}^N \in \Delta^N$ and then considering the following scenario-based NE seeking problem, where each agent $m \in \mathcal{M}$ solves the following optimization program
\begin{align}
 \min_{x^m \in X^m}{J_m(x^m,x^{-m})} \ \text{subject to} \  x^m \in \bigcap_{i=1,\dots,N} X^m _{\delta_i}. \label{EVmeta}
\end{align}
Our aim is to provide probabilistic feasibility guarantees for the entire set of NE of (\ref{EVmeta}) returned by an arbitrary NE seeking algorithm, i.e., to quantify the probability that any NE strategy $x^m_{NE}$ of (\ref{EVmeta}) belongs to the constraint set $X^m_\delta, \ \forall \ m\in \mathcal{M}$ for a new unseen sample $\delta \in \Delta$. To this end, our analysis is primarily focused on feasibility problems affected by uncertainty. Based on the derived results we revisit problem (\ref{EVprin}) and attempt to provide robustness certificates for the obtained NE for a motivating application which fits in this class of games, i.e, the electric vehicle (EV) charging control problem, presented below. 
\subsection{EV-charging control problem}
The EV-charging control problem can be treated as a noncooperative game comprised of self-interested agents-vehicles each of them aiming at minimizing their own electricity cost, while their charging schedules are subject to certain specifications. The  two main requirements for the operation of the system, namely, the lower and upper bounds imposed on the charging schedule and the total energy level to be achieved at the end of charging, can be modeled as constraints of affine form. However, most of the work up to this point assumed that these constraints are purely deterministic \citep{Callaway2013}, \citep{Paccagnan2019b}, \citep{Deori2018}. We extend this framework by imposing uncertainty on the constraints. In this case, each agent $m \in \mathcal{M}$ solves the following problem 
\begin{align}
&\min_{x^m \in \mathbb{R}^{d}} J_m(x^m,x^{-m})  \ \text{subject to}  \nonumber \\
&x^m \in \bigcap\limits_{\delta \in \Delta}\{ [\underline{x}^m (\delta),\overline{x}^m(\delta)]\bigcap\{x^m \in \mathbb{R}^{d}: \sum_{t=1}^{d}x_t^m \geq E^m(\delta) \} \} .  \label{evgames}
\end{align} 
The variables $x^m=(x^m_t)_{t=1}^d$ and $J_m$ denote, respectively, the charging schedule for all time instances $t \in \{1, \dots, d\}$ and the electricity cost to be minimized for each vehicle $m \in \mathcal{M}$. The uncertain constraint for each vehicle $m$ comprises of uncertain lower and upper bounds  $\underline{x}^m (\delta)$ and $\overline{x}^m(\delta)$ and total energy levels  $E^m (\delta)$.
The uncertainty inherent in those constraints stems from several factors such as the battery dynamics of each vehicle, the preferences of the users and the status of the electric grid to name a few. Due to the presence of a variety of unpredictable internal and external influences contributing to the uncertainty of the system, it is very difficult to address the problem using traditional probabilistic approaches. We thus adopt a data-based approach, a more viable alternative. \par

\section{Probabilistic guarantees of feasibility programs}
Motivated by the EV-charging control problem discussed in the previous section, we leverage the recent results of the scenario approach \citep{CampiGarattiRamponi2018} to provide probabilistic guarantees for the (possibly) multiple NE of (\ref{evgames}) in a computationally efficient manner. As NE constitute feasible solutions as far as constraint satisfaction is concerned, we focus on obtaining robustness certificates for a more general class of problems, that of feasibility problems under uncertain constraints, where the NE problem (\ref{evgames}) emerges as a special case. As such, we consider the following feasibility problem:
\begin{align}
\mathrm{P}_{\Delta}: \mathrm{find} \ {x \in X}, \ 
&\text{subject to } x \in \bigcap_{\delta \in \Delta} X_\delta, \nonumber
\end{align}
where $x$ is a decision vector belonging to the set $X \subset \mathbb{R}^d$ and $\delta$ is a random variable, defined as in Section 2.1,  that encodes the uncertainty parameterizing agents' constraint sets. Note that for the gaming set-up, as introduced in Section 2, $X= \mathbb{X}=\prod_{m \in \mathcal{M}} X^m$ and $X_\delta=\prod_{m \in \mathcal{M}} X^m_{\delta}$.
For the feasibility problem under study, we define the following scenario program:
\begin{align*}
\mathrm{P}_{N}: \mathrm{find} \ {x \in X} \
&\text{subject to } x \in \bigcap_{i=1,...,N} X_{\delta_i},
\end{align*}
where $N$ denotes the number of samples $\delta_i, i=1,..., N$, drawn in an i.i.d. fashion from $\Delta$ according to the probability distribution $\mathbb{P}$ and  $X_{\delta_i}$ their respective constraints. 
Our results depend on an affine constraint structure, thus we impose the following standing assumption: 
\begin{assumption} \label{ass1} \hfill{}\\[-3ex]
\begin{enumerate}
\item The deterministic constraint set $X$  is a non-empty, compact and convex polytope\footnote{A polytope $\Pi \in  \mathbb{R}^d$ can be expressed by its H-representation, i.e., the intersection of a finite number of halfspaces, and also as the convex hull of its vertex set $v(\Pi)=\{x_1,...,x_Q \}  $ i.e, $\Pi=conv(v(\Pi))=\{ \sum_{j=1}^{Q}{x_j\lambda_j} : \sum_{j=1}^{Q}{\lambda_j}=1,  \lambda_j \geq 0, \ j=1,...,Q\}$, where $v(\cdot)$ and $conv(\cdot)$ denote the set of vertices of the polytope  and the convex hull, respectively. This representation is generally known as $V$-representation.}.
\item  For each $\delta \in \Delta$, $X_{\delta}=\{ x \in \mathbb{R}^d : g(x,\delta) = a^Tx - b \leq 0 \}$, where $g$ is an affine function given by the mapping $g: X \times \Delta \rightarrow \mathbb{R}$, where $a \in \mathbb{R}^{d}$, $b \in \mathbb{R}$ and $\delta=(a^T \  b) \in \mathbb{R}^{d+1}$.
\item For each multi-sample $\{\delta_i\}_{i=1}^N$ the polytope \\ $\Pi_N= \{ \bigcap_{i=1}^N  X_{\delta_i} \} \bigcap X= \{x \in X :g(x,\delta_i) \leq 0, i=1,..., N \}$  has a non-empty interior.
\end{enumerate}
\end{assumption}
Note that vector valued affine functions are also captured by our framework; see example of Section 4.
Assumption \ref{ass1} guarantees that  the polytope $\Pi_N$  is compact and $P_N$ admits at least one solution for any chosen multisample  $\{\delta_i\}_{i=1}^N$. Under Assumption \ref{ass1} the feasibility problem $P_N$ can be equivalently written as 
\begin{align*}
\mathrm{P}_{N}: \mathrm{find} \ x \ 
&\text{subject to } x \in \Pi_N.
\end{align*}
Upon finding the feasibility domain $\Pi_N$ of the problem $\mathrm{P}_{N}$, we are interested in investigating the robustness properties collectively for all the points of this domain to yet unseen samples, in other words in quantifying the probability that a new sample $\delta \in \Delta$ is drawn such that the constraint $X_\delta$ defined by this sample is not satisfied by some given point $x\in \Pi_N$. This concept, which is of crucial importance for our work, is known in the literature as the probability of violation and  is adapted in our context to represent the probability of violation of a set. By Definition 1 in \citep{CampiCalafiore2006} the probability of violation of a given point $x \in \Pi_N$ is defined as
\begin{align}
& V(x)=  \mathbb{P} \Big \{ \delta \in \Delta :~ x \notin X_\delta\Big \}. \label{Violationpoint}
\end{align} 
By Assumption \ref{ass1}, the probability of violation can be equivalently written as $V(x)=  \mathbb{P} \Big \{ \delta \in \Delta :~ g(x,\delta)>0 \Big \}$. We can now define the probability of violation of the polytope $\Pi_N$.
\begin{definition}  \label{def2}
Let  $\mathcal{P} \subseteq 2^X$  be the set of all non-empty, compact and convex polytopes that are subsets of $X$. For any $\Pi_N \in \mathcal{P}$ we define  the probability of violation of the set $\Pi_N$  as  a mapping $\mathbb{V}: \mathcal{P} \rightarrow [0,1]$ given by the following relation: 
\begin{align*} 
\mathbb{V}(\Pi_N)&= \sup_{x \in \Pi_N} V(x).
\end{align*} 
\end{definition} 
Definition \ref{def2} can be considered a special case of Definition 2 in  \citep{GrammaticoZMGL14}. In Definition \ref{def3} three concepts of crucial importance are introduced.
\begin{definition} \label{def3}
\begin{enumerate}
\item  For any $N$, an algorithm is a mapping $ A_N: \Delta^N \rightarrow \mathcal{P}$ that associates the multisample  $\{\delta_i\}_{i=1}^N$ to a unique polytope $\Pi_N \in \mathcal{P}$. 
\item  Given a multisample $\{\delta_i\}_{i=1}^N \in \Delta^N$, a support subsample $S \subseteq \{\delta_i\}_{i=1}^N$ is a subset of the entire multisample with cardinality $k \leq N$ so that for $ i_1<i_2<...<i_k$ and $I_k=\{i_1,i_2,...,  i_k\}$,  $S=\{\delta_i\}_{i \in I_k}$ is such that  $A_k(\{\delta_i\}_{i \in I_k})=A_N(\{\delta_i\}_{i=1}^N)$, i.e., the solution returned by an algorithm when fed with the subsample is the same with the one obtained when the entire multisample is used.
\item  A support subsample function is a function of the form $B_N: \{\delta_i\}_{i=1}^N  \rightarrow \{ i_1,..., i_k\}$  that takes as input all the samples and returns as output the indices of only these samples  that constitute an irreducible\footnote{A support subsample $S=\{\delta_i\}_{i \in I_k} \subseteq \{\delta_i\}_{i=1}^N$  is said to be irreducible if no element can be further removed from S leaving the solution unchanged.} \label{footnote1} support subsample. \\
\end{enumerate}
\end{definition}
 Note that the notions of support subsample and support subsample function in Definitions \ref{def3}.2, \ref{def3}.3 are respectively referred to as compression set and compression function in \citep{Margellos2015}.  Moreover, Definition \ref{def3}.3 in our case translates into the relation $ \Pi_N=\Pi_{\{\delta_i\}_{i \in I_k}}$, where the cardinality of the support subsample $\{\delta_i\}_{i \in I_k}$ is by definition the cardinality of the output of the function $B_N$. \par 
Let $K_N=|v(\Pi_N)|$ and $F_N$ be the number of vertices and the number of 
facets\footnote{For a definition of the facets of a polytope we refer the reader to Definition \ref{facets} of the Appendix.} of $\Pi_N$, respectively. As shown in Theorem 8.2(b) of \citep{Bronsted1982}, $F_N$ is finite. It is important to emphasize that the dependence of the polytope $\Pi_N$ on the multi-sample $\{\delta_i\}_{i=1}^N$ implies that both  $F_N$ and $K_N$ are random variables that depend on $\{\delta_i\}_{i=1}^N$. 

Next we define the set 
\begin{align}
 \mathcal{P}_\delta&=\{\Pi \in \mathcal{P}:  g(x,\delta) \leq 0, \  \forall \ x \in v(\Pi) \} \nonumber \\
 &=\{\Pi \in \mathcal{P} : \Pi \subseteq X_\delta \} \label{P},
\end{align}
of all the non-empty, compact and convex polytopes $\Pi$ that satisfy the constraint associated with the sample $\delta \in \Delta$. Note that if all the vertices of the polytope satisfy the inequality $g(\cdot,\delta) \leq 0$, then every point $x\in \Pi$ of the polytope satisfies it as well, since $x$ can always be expressed as a convex combination of the polytope's vertices. \par 

 We are now ready to introduce the following theorem, which is the main result of our paper.

\begin{thm} \label{Theorem1}
Consider Assumption \ref{ass1} and any $A_N, B_N$ as in Definition \ref{def3}. Fix $\beta \in (0,1)$ and define the violation level  $\epsilon: \{0,...,N\} \rightarrow [0,1]$ as a function such that
\begin{align}
 \epsilon(N)=1 \ \text{and} \  \sum_{k=0}^{N-1} {N \choose k}(1-\epsilon(k))^{N-k}=\beta. \label{epsilon}
\end{align} 
 We have that 
\begin{align*}
\mathbb{P}^{N} \Big \{\{\delta_i\}_{i=1}^N \in \Delta^{N}:~\mathbb{V}(\Pi_N)> \epsilon(F_N) \Big \} \leq \beta,
\end{align*}
 where $\mathbb{P}^{N}= \prod_{i=1}^N \mathbb{P}$  is the product probability measure, and $F_N$ is the number of facets of $\Pi_N$.
\end{thm}

\emph{Proof}:
The first part of the proof closely follows that of Theorem 1 in \citep{GrammaticoZMGL14}. For a fixed multisample  $ \{\delta_i\}_{i=1}^N \in \Delta^{N}$ consider an arbitrary point $x \in \Pi_N$. Then,  the following inequalities are satisfied
\begin{align}
&  V(x)  =  \mathbb{P} \Big \{ \delta \in \Delta :~ x  \notin X_\delta \Big \}=  \mathbb{P} \Big \{ \delta \in \Delta :~g( x,\delta)>0\Big \} \nonumber \\
&\stackrel{(i)}{=} \mathbb{P} \Big \{ \delta \in \Delta :~ g( \sum_{j \in I_{d+1}} {\lambda_j x_j},\delta)>0\Big \}  \nonumber \\
& \stackrel{(ii)}{=} \mathbb{P} \Big \{ \delta \in \Delta :~ \sum_{j \in I_{d+1}}\lambda_j  g( x_j,\delta)>0\Big \}  \nonumber \\ 
& \leq \mathbb{P} \Big \{ \delta \in \Delta :~  \sum_{j \in I_{d+1}}\lambda_j \max_{j \in I_{d+1}} g( x_j,\delta)>0\Big \}  \nonumber \\
& \leq \mathbb{P} \Big \{ \delta \in \Delta :~ \max_{j \in I_{d+1}} g( x_j,\delta) >0\Big \}   \nonumber  \\
& = \mathbb{P}\bigg \{ \bigcup_{j \in I_{d+1}} \Big \{ \delta \in \Delta :~ g( x_j,\delta) >0\Big \} \bigg \}  \nonumber   \\
&\stackrel{(iii)}{\leq} \mathbb{P}\bigg \{ \bigcup_{j=1}^{K_N} \Big \{ \delta \in \Delta :~ g( x_j,\delta) >0\Big \} \bigg \}.  \label{result1}
\end{align} 
 Equality (i) is derived from Caratheodory's Theorem (Theorem 1.1 in \citep{Monroy2015}) where the set under study is the polytope $\Pi_N$. In our case, Caratheodory's Theorem states that any arbitrary point of the polytope $x \in \Pi_N$  can be represented as a convex combination of at most $d+1$  vertices from the set $v(\Pi_N)$, which means that there exists a subset of indices  $ I_{d+1} \subseteq \{1,..., K_N\} $  such that  $x=\sum_{j \in I_{d+1}} {\lambda_j x_j}$, where $\sum_{j \in I_{d+1}} \lambda_j=1$ and $\lambda_j \geq 0, \ \forall \ j \in I_{d+1}$.  Equality (ii) stems from the fact that g is an affine function of $x$ for any given $\delta \in \Delta$ due to Assumption \ref{ass1}. The last inequality follows from the fact that $I_{d+1} \subseteq \{1,...,K_N\}$, due to the fact that $K_N \geq d+1$ as the polytope has a non-empty interior by Assumption 3 (3). Since (\ref{result1}) holds for all $x \in \Pi_N $, we have that 
\begin{equation}
 \mathbb{V}(\Pi_N)=\sup_{x \in \Pi_N} V(x) \leq \mathbb{P}\bigg \{ \bigcup_{j=1}^{K_N} \Big \{ \delta \in \Delta :~ g( x_j,\delta) >0\Big \} \bigg \}. \nonumber
\end{equation}
Therefore, for any multisample $\{\delta_i\}_{i=1}^N$ and for any cardinality (not necessarily irreducible) of the support subsample $k \in \{1,...,N\}$ the following inequalities are satisfied:
\begin{align} 
\mathbb{P}^{N} \Big \{  \{\delta_i\}_{i=1}^N \in \Delta^{N} &:~ \mathbb{V}(\Pi_N) > \epsilon(k) \Big \} \nonumber \\
\leq \mathbb{P}^{N} \Big \{  \{\delta_i\}_{i=1}^N \in \Delta^{N}& \!\!:~ \nonumber \\
 \!\!\!\! \mathbb{P}\Big \{ \bigcup_{j=1}^{K_N}  \Big \{& \delta \in \Delta : g(x_j,\delta)>0 \Big \}\Big \}> \epsilon(k) \Big \} \nonumber\\
= \mathbb{P}^{N} \! \Big \{\! \{\delta_i\}_{i=1}^N\! \!\in \Delta^{N}\! :\!\!\!\!~& \nonumber \\
 \mathbb{P} \Big \{ \!\!\delta \in \Delta\! : \ & \exists \ \!\!  x \in v(\Pi_N)\!\!, g(x,\delta)\!\!>\!\! 0 \Big \}\!\!\!>\!\epsilon(k) \Big \} \nonumber\\
=\mathbb{P}^{N} \Big \{  \{\delta_i\}_{i=1}^N \in \Delta^{N}&:~ \nonumber \\
\mathbb{P} \Big \{ \delta \in \Delta& :~ \Pi_N \not\subseteq X_\delta\Big \}  > \epsilon(k) \Big \}, \label{first}
\end{align}
where the last inequality is due to (\ref{P}). Define now an algorithm $A_N$ as in Definition  \ref{def3}.1, that returns the polytope confined by the feasibility region of $\mathrm{P}_{N}$. By construction, $A_N$ satisfies Assumption 1 of \citep{CampiGarattiRamponi2018},  since for any multisample $\{\delta_i\}_{i=1}^N$ it holds that $A(\{\delta_i\}_{i=1}^N) \in \mathcal{P}_{\delta_i}$, for all $i=1,..,N$. The satisfaction of Assumption 1  paves the way for the use of Theorem 1 of \citep{CampiGarattiRamponi2018}. 
In particular, Theorem 1 of \citep{CampiGarattiRamponi2018} implies that the right-hand side of (\ref{first}) can be upper bounded by $\beta$, for $k$ being the cardinality of an irreducible support subsample (see Definition \ref{def3}.1). For our case, an irreducible (in fact minimal) subsample coincides with the minimum number of facets that construct the polytope $\Pi_N$, i.e., $k=F_N=\mathrm{rank}([H \  L])$, where matrices $H$, $L$ are of appropriate dimension and constitute the H-representation of $\Pi_N$, i.e., $\Pi_N=\{x \in \mathbb{R}^{d}: Hx \leq L\}$. As such,  we have that
\begin{align}
&\mathbb{P}^{N} \Big \{ \{\delta_i\}_{i=1}^N \in \Delta^{N}\!\!:\!\!~\mathbb{P} \Big \{ \delta \in \Delta :~\!\!\!\! \Pi_N \not\subseteq X_\delta\Big \}  > \epsilon(F_N) \Big \} = \nonumber \\
&\mathbb{P}^{N} \Big \{ \{\delta_i\}_{i=1}^N \in \Delta^{N}\!\!:\!\!~\mathbb{P} \Big \{ \delta \in \Delta :~\!\!\!\! \Pi_N \notin \mathcal{P}_\delta \Big \}  > \epsilon(F_N) \Big \} \leq \beta. \label{second}
\end{align}
From (\ref{first}) and (\ref{second}) we obtain that:
\begin{align}
 \mathbb{P}^{N} \Big \{  \{\delta_i\}_{i=1}^N \in \Delta^{N}:~\mathbb{V}(\Pi_N) > \epsilon(F_N) \Big \}  \leq \beta,
\end{align}
 thus concluding the proof.
 \qed \\
Note that, even though within the proof of our theorem we also use the V-representation of the polytope, only the number of facets are needed  to provide probabilistic guarantees for the entire feasibility region. This feature is appealing from a computational point of view as, in most practical cases, the constructed polytope has a significantly smaller number of facets than vertices.
To illustrate this, consider the EV charging control problem of Section 2. Vehicles' charging schedule is subject to upper and lower bounds at each time instance $t \in \{1,\ldots,d\}$. Hence, for a multi-sample $\{\delta_i\}_{i=1}^N \in \Delta^N$,  the feasibility domain $\prod \limits_{ t \in \{ 1,\dots,d \}}\prod \limits_{m \in \mathcal{M}}\bigcap\limits_{i=1,\dots,N}[\underline{x}_t^m (\delta_i),\overline{x}_t^m(\delta_i)]$ of the problem is a hyperrectangle whose number of facets $F_N=2Md$ grows linearly with respect to the number of decision variables, while the number of vertices is given by $K_N=2^{Md}$, which grows at an exponential rate with respect to $d$. In the Appendix we provide a relationship between $\mathbb{V}(\Pi_N)$ and the probability of constraint violation of the polytope's vertices. This is not used further in the paper but is interesting per se.  \par
Quantifying the cardinality of the minimal support subsample becomes trivial, since in our case it coincides with the number of facets, thus circumventing the need of employing the greedy algorithm in \citep{CampiGarattiRamponi2018}. 
A direct consequence of Theorem \ref{Theorem1} is that we can provide distribution-free guarantees for any NE of problem (\ref{EVprin}). To solidify this statement we introduce the following corollary
\begin{cor}
Consider Assumption \ref{ass1} and the setting of Theorem \ref{Theorem1}.
 We have that 
\begin{align*}
\mathbb{P}^{N} \Big \{\{\delta_i\}_{i=1}^N \in \Delta^{N}:~ V(x_{NE})> \epsilon(F_N), \  \text{for any} \  x_{NE}  \ \text{of (\ref{EVmeta}}) \Big \} \leq \beta.
\end{align*}
\end{cor}
Note that the choice of algorithm to determine a NE of the problem is arbitrary.

\section{Numerical examples}
\subsection{Feasibility of a random 2-dimensional polytope}

We initially apply our results to a 2-dimensional example of a polytope constructed by the intersection of random halfspaces of the form $a_1x_1+a_2x_2-b \leq 0$, where $a_1, a_2$ and $b$ are scalars following uniform distributions with support $[-4, 4]$, $[-4, 4]$ and $[10, 15]$, respectively. Each sample $\delta$ is defined as a vector $\delta=(a_1, a_2,  b) \in \mathbb{R}^{3}$. \par 
The theoretical relation between the violation level $\epsilon(k)$ for different values of $k$ is illustrated in Figure 1, where $\epsilon$ is computed according to (\ref{epsilon}) by fixing the confidence parameter to $\beta=10^{-6}$. The colour code corresponds to different choices of the number of samples $N$.  As it can easily be observed, choosing a larger number of samples improves the robustness guarantees of the feasibility region significantly, as lower $\epsilon(k)$ implies lower probability of constraint violation. \par

\begin{figure} 
\hspace*{-0.3cm}  
\includegraphics[width=9.5cm, height=6cm, trim={0cm, 1, 0cm 0},clip]{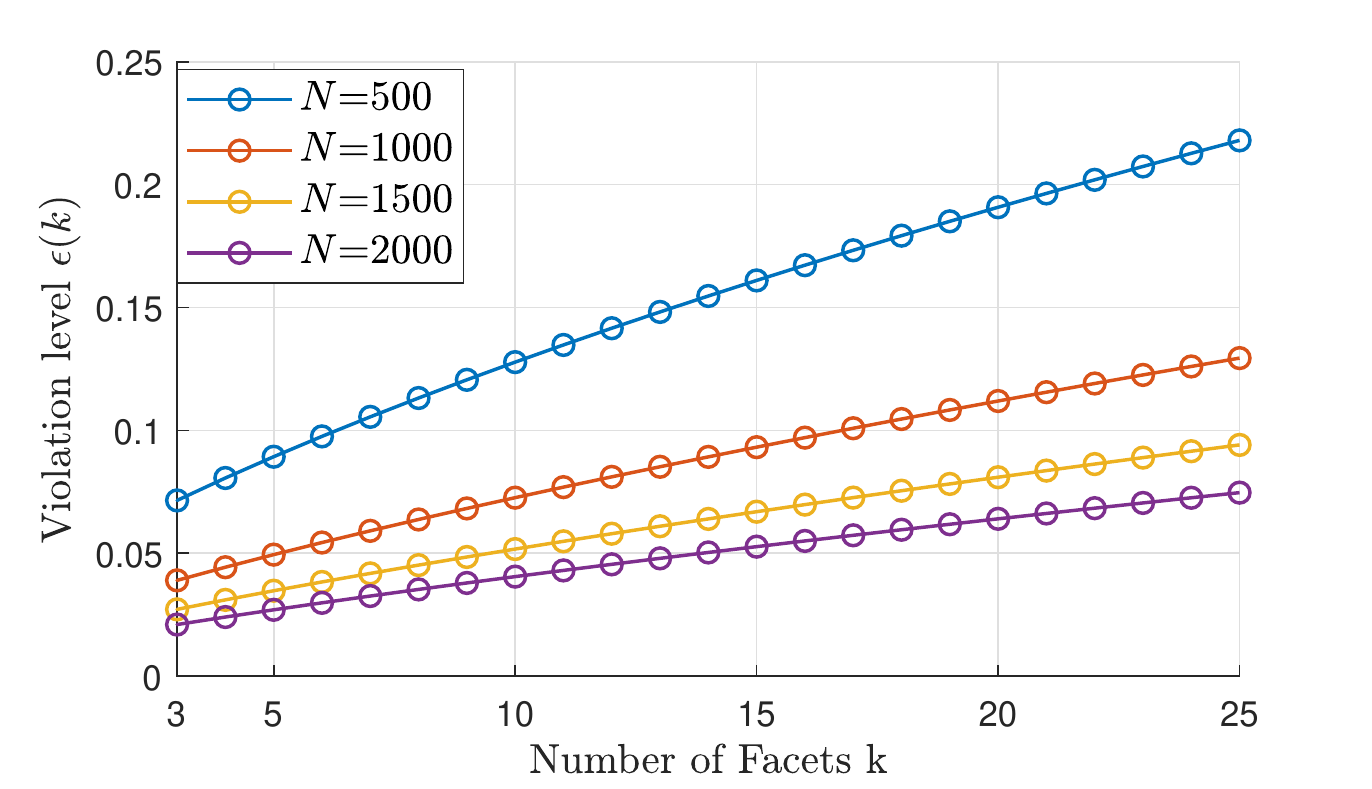}
\caption{ \small The violation level $\epsilon(k)$ as a function of the empirically most frequently observed numbers of facets $k$. Confidence level is set to $\beta = 10^{-6}$, while four different choices for the number of samples $N$ are investigated.}
\end{figure}

To test the validity of our theoretical bounds in practice we need to compute the probability of violation $\mathbb{V}(\Pi_N)$ and compare it with the guarantees provided by Theorem \ref{Theorem1}.  Let $\mathbb{\hat{V}}(\Pi_N)$ denote an empirical estimate of $\mathbb{V}(\Pi_N)$ and $\hat{V}(x)$ an empirical estimate of a point $x \in \Pi_N$. By gridding the polytope using a large enough number of points $x_r, r=1,...,R$ that cover the entire polytope we have that  $\mathbb{\hat{V}}(\Pi_N) =\max_{r=1,...,R} V(x_r)$.  A sufficiently large number of test samples is used, different from those used to construct the polytope. \par

\begin{figure}
\includegraphics[width=10cm, height=7cm]{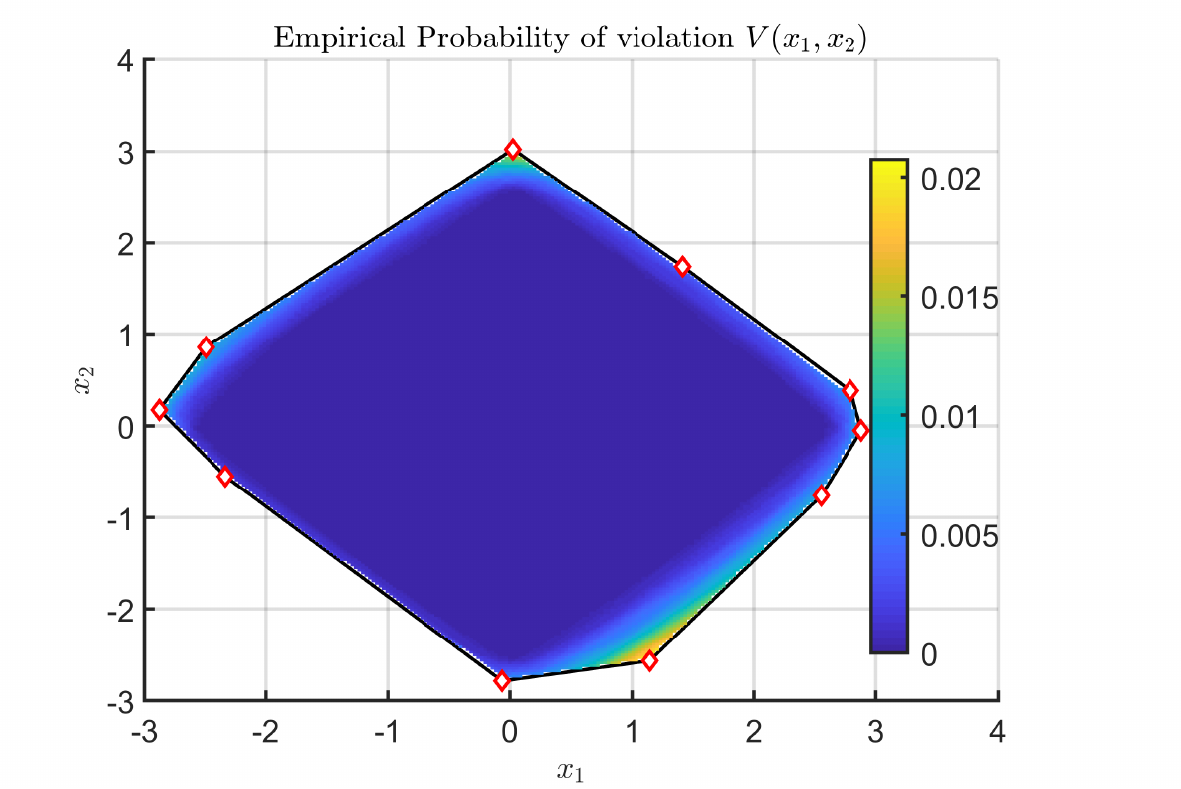}
\caption{\small  The empirical probability of violation of each point of the polytope under study constructed by 100 random realizations of the affine constraints $a_1x_1+a_2x_2-b \leq 0$. The number of grid points is $R=48633$, while $10000$ test samples are used.  Note that the higher probability of violation occurs at one of the vertices (red diamonds).}
\end{figure} 

To this end, we generate a total number of $N=100$ samples for the construction of the random polytope shown in Figure 2 and use 10000 test samples to compute the empirical probability of violation for each point of the grid.  The highest probability of violation occurs at a polytope vertex; Lemma 1 provides some theoretical support to this numerical evidence, by showing that the worst-case probability of violation is proportional to the worst-case vertex violation. A tighter relation for certain class of programs is currently under investigation. \par

Finally, we validate the derived theoretical bounds against the empirical probability of violation of 50 independent realizations of polytopes. This means that each polytope is constructed using a different multi-sample $\{\delta_i \}_{i=1}^{2000}$. By keeping the same value for $\beta$, as in the previous case, we count the number of facets $F_N$ of each polytope and then compute the theoretical bound of the violation level that corresponds to it. Subsequently, using 20000 test samples we compute an empirical estimate of the probability of violation for each polytopic realisation, as outlined above. If there is more than one polytope among the 50 that has the same number of facets, we choose the one with the maximum empirical probability of violation. As anticipated, $\epsilon(k)$ constitutes an upper bound for any of the computed empirical probabilities of violation.

\begin{figure}
\hspace*{-0.3cm}  
\includegraphics[width=8.7cm, height=6cm, trim={0cm, 3, 1cm, 1},clip]{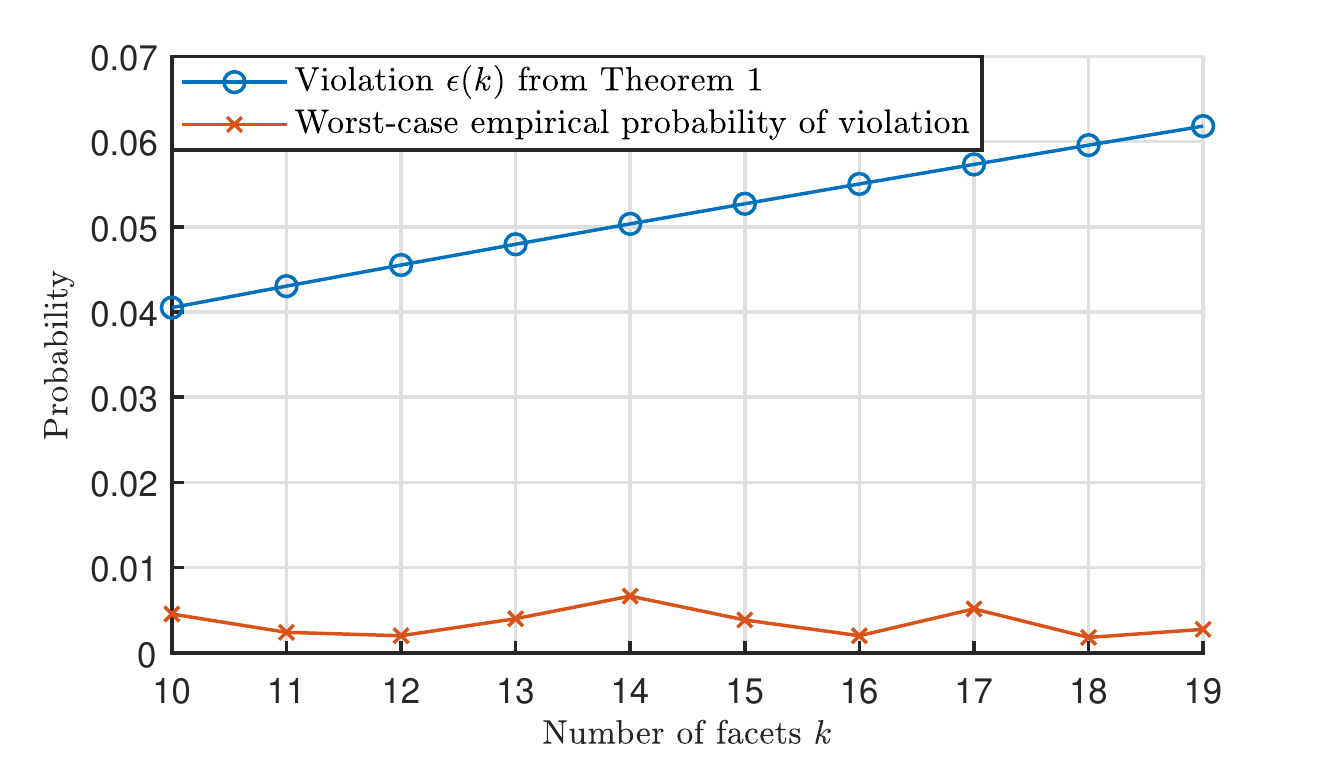}
\caption{\small Comparison of the theoretical $\epsilon(k)$ derived from Theorem \ref{Theorem1} (blue) and the worst-case empirical probability of violation of the entire polytope that corresponds to each $k$ (red). We generated 50 different polytopes, each of them using 2000 samples, while for each one the worst case probability corresponds to the grid point with the highest probability of constraint violation. The confidence level was fixed to $\beta=10^{-6}$. Note that the non-monotonic behaviour for the empirical probability of violation is due to the fact that for any number of facets $k$ the number of polytopes among which the worst-case probability of violation is calculated is not the same. }
\end{figure}

\subsection{PEV charging game revisited}
We revisit the EV charging control game of Section 2.2. Our aim is to provide guarantees on the probability that a NE satisfies the constraints of (3). We assume that the upper constraint  $(\overline{x}^m(\delta_u))_{m\in \mathcal{M}}\in \mathbb{R}^{Md}$ is affected by an additive uncertainty in the form of  $\delta_{u} \in \mathbb{R}^{Md}$, whose elements are random realisations of $0.3U(0,1)\cdot \mathcal{N}(1,3)$, where  ${U(0,1)}$ is a random variable that follows a uniform distribution with support $[0, 1]$ and $\mathcal{N}(1,3)$ another random variable that follows the gaussian distribution with mean $1$ and standard deviation $3$. As such, $\overline{x}^m(\delta_u)=\overline{x}^{nom}+\delta_{u}$, where $\overline{x}^{nom}$ is a given deterministic component. The total energy $E=(E^m)_{m \in \mathcal{M}}$ of each agent at the end of charging is also affected by uncertainty i.e., $E^m=(1-\delta_e)E^m_{nom}$, where $\delta_e  \in \mathbb{R}^{M}$ and its elements are extracted according to the probability $0.05\mathcal{N}(0,1)$ and $E^m_{nom} \in \mathbb{R}$ is the nominal final energy demand of each agent $m \in\mathcal{M}$ drawn from $U(10,17)$. The uncertainty vector is given by $\delta=[\delta_{u}, \delta_{e}]\in \mathbb{R}^{M(d+1)}$. The lower bound is assumed to be deterministic and, particularly, $\underline{x}^i=0$ for any $ i\in \mathcal{M}$. Finally, the cost function of each vehicle $m \in \mathcal{M}$ is given by $J_m(x^m,x^{-m})=(x^m)^T(A_0\sigma(x^m,x^{-m})+b_0)$, where the matrix $A_0 \in \mathbb{R}^{d\times d}$ is diagonal and $\sigma(x^m, x^{-m})=\sum_{m=1}^M x^m$. Following the work of \citep{Fele2019a}, \citep{Fele2019b} the entries $\{a_t\}_{t =1 }^d$ of $A_0=\text{diag}(\{a_t\}_{t =1 }^d)$ are evaluated by rescaling a winter weekday demand profile in the UK \citep{NationalGrid}. The vector $b_0 \in \mathbb{R}^d$ is set to zero. \par

\begin{figure} [t]
\hspace*{-0.29cm}  
\includegraphics[width=9.3cm, height=7cm, trim={0cm, 1, 0cm, 1},clip]{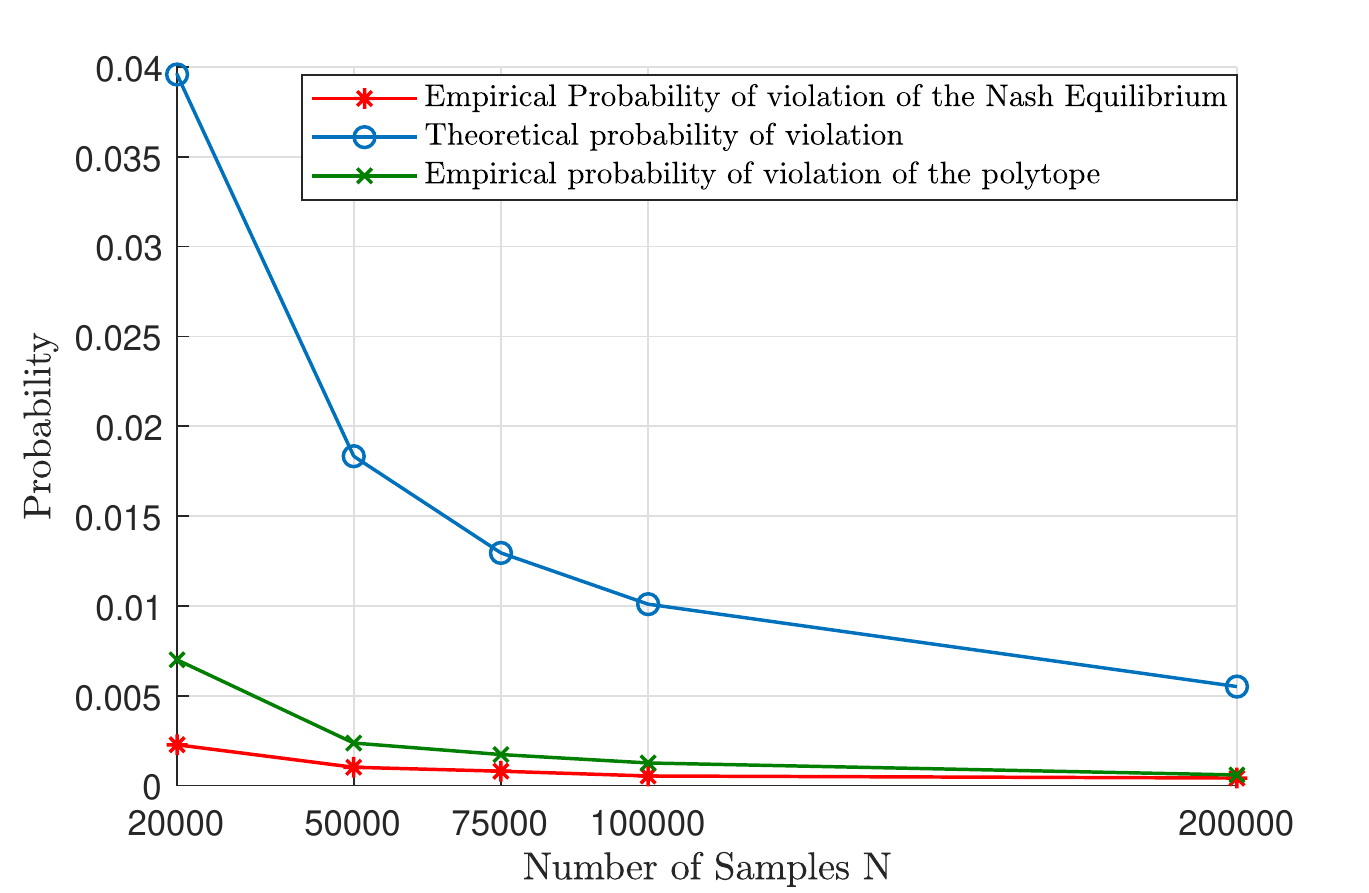}
\caption{\small The empirical probability of violation of the entire polytope $\mathbb{\hat{V}}(\Pi_N)$ (green) and the empirical probability of violation of the computed NE returned by the algorithm of \citep{Fele2019a},\citep{Fele2019b}, namely, $\hat{V}({x_{NE}) }$  (red) versus the theoretical violation level of Theorem \ref{Theorem1} (blue) with respect to five different values of the number of samples $N=20000, 50000, 75000, 100000, 200000$. The empirical probability of violation for both the polytope and the NE is computed using 2000000 test samples. Note that the blue line corresponds to the theoretical counterpart of the green one.}
\end{figure}

\begin{figure} [!htbp]
\hspace*{-0.3cm}  
\includegraphics[width=9.3cm, height=7cm, trim={0cm, 1, 0cm, 1},clip]{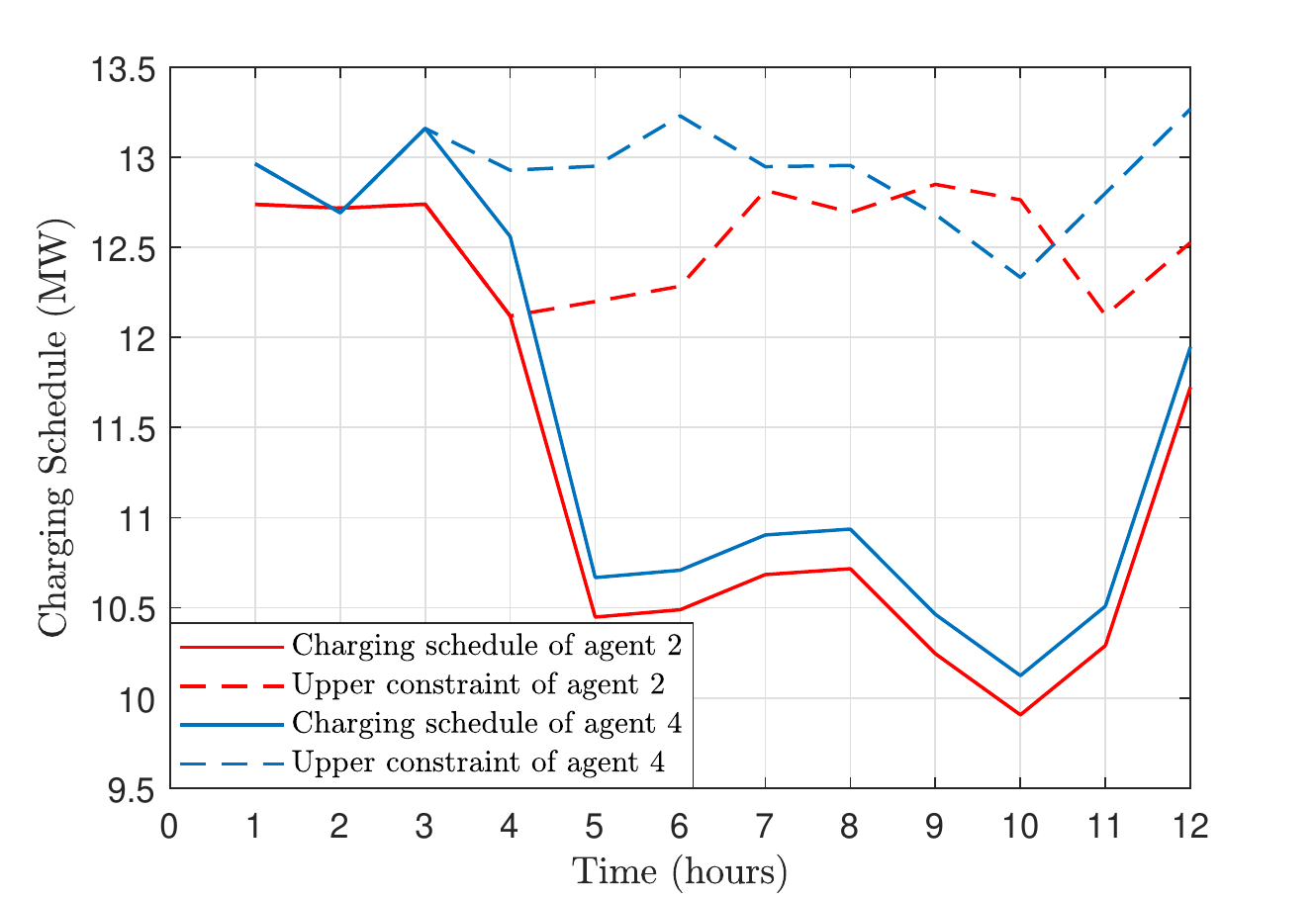}
\caption{\small Charging schedules of agents 2 (red solid line) and 4 (blue solid line) and their respective uncertain upper bounds. The agents' schedules and their upper constraints coincide during the beginning and the end of the horizon, which implies that they are charging at their boundaries in their attempt to benefit from the price die during these time instances.}
\end{figure}

To compute the NE of (3) we employed the algorithm of \citep{Fele2019b} (see Algorithm 1 therein). Setting the number of agents to $M=10$, the number of timesteps to $d=12$ and $\beta=10^{-5}$, we run the algorithm for multi-samples of different size, namely, 20000, 50000, 75000, 100000 and 200000 and we compare the behaviour of two different probabilities  with respect to the cardinality of the multisample, as illustrated in Figure 4.  With green we illustrate the worst case empirical probability of violation for the polytope $\hat{\mathbb{V}}(\Pi_N)$ (calculated by taking the maximum violation among all grid points on the polytope), while with red we show the empirical probability $\hat{V}(x_{NE})$ of violation for the $x_{NE}$ returned by the algorithm of \citep{Fele2019a},\citep{Fele2019b}.  Note that the blue line is the theoretical counterpart of the green line and corresponds to the theoretical violation level $\epsilon(k)$, as defined in Theorem 6. The empirical calculation was performed using 2000000 test samples, different from those used in the NE seeking process. As expected, both empirical values are less than the theoretical bound derived by Theorem \ref{Theorem1}. It was also anticipated that for any multisample we would have that  $\mathbb{\hat{V}}(\Pi_N) \geq \hat{V}(x_{NE})$ as the former corresponds to the collective violation of all feasible points, including $x_{NE}$. Finally, Figure 5 illustrates the charging schedules of two of the agents and their respective uncertain upper bounds. Due to the high required total energy at the end of charging, we observe that the agents' schedules and their upper constraints coincide during the beginning and the end of the horizon, which implies that they are charging at their boundaries in their attempt to exploit the low price during these time instances and minimize their charging cost.
\section{Concluding remarks}
Considering a feasibility problem under uncertain polytopic constraints, we provided probabilistic guarantees for the entire feasibility set in an \emph{a posteriori} fashion. The importance of this result is better shown in the context of the EV-charging control problem where computationally efficient robustness certificates are obtained for the NE returned by any algorithm. Effort is being made towards extending our results to include feasibility problems subject to uncertain convex constraints and in investigating the case  of uncertain games where each agent's samples are drawn from her own private uncertainty set, rather than from a common set of samples. Finally, simulation results indicate the validity of a stronger statement than that of Lemma \ref{lem1}, that is, the highest probability of violation among all the points of the polytope occurs at (at least) one of its vertices for some problems; we aim to investigate theoretically the class of programs for which this indication is valid.
\section{Appendix} 
We provide some auxiliary definitions for the analysis of Section 3. 
\begin{definition} (valid inequalities, faces and facets) \citep{Ziegler1995} \label{facets}
\begin{enumerate}
\item Let $\Pi \subset \mathbb{R}^{d}$ be a convex polytope. An affine inequality $a^Tx \leq b$, where $a$, $b$ are of appropriate dimensions, is valid for $\Pi$, if it is satisfied for all points $x \in \Pi$. 
\item A face of $\Pi$ is  defined as any set of the form 
\begin{align}
&f = \Pi \cap \{x \in \mathbb{R}^{d} : a^Tx = b\}, \ \text{where} \nonumber
\\ & a^Tx \leq b \ \text{is a valid inequality for} \  \Pi. \label{face}
\end{align}
\item If alongside (\ref{face}), the additional condition $\text{dim}(f)= \text{dim}(\Pi)-1$ is satisfied, where $\text{dim}(f)$, $\text{dim}(\Pi)$ denote the dimensions of the face $f$ and the polytope $\Pi$, respectively, then the face is also referred to as a facet of the polytope. 
\end{enumerate}
\end{definition}
The following lemma provides an additional result, which even though it is not used for our derivations, is interesting per se. 
\begin{lemma} \label{lem1}
Consider the sets $X$ and $X_{\delta_i},  i=1,..., N$, that satisfy Assumption \ref{ass1}. Then for any given multisample $ (\delta_1,...,\delta_N) \in \Delta^{N}$, there exists a vertex $\hat{x} \in v(\Pi_N)$ such that  that: 
\begin{align*}
& \mathbb{V}(\Pi_N) \leq (d+1) \mathbb{P} \Big \{ \delta \in \Delta :~ \hat{x} \notin X_\delta\Big \}=(d+1)V(\hat{x}).\\
\end{align*}
\end{lemma}
\emph{Proof}:
Consider a fixed multisample and any arbitrary point $x \in \Pi_N$. By relation (iii) of (7)  we can follow an alternative direction which leads to the following inequalities:
\begin{align*}
&V(x) \leq  \mathbb{P}\bigg \{ \bigcup_{j \in I_{d+1}} \Big \{ \delta \in \Delta :~ g( x_j,\delta) >0\Big \} \bigg \}  \\
&\leq \sum_{j \in I_{d+1}}\mathbb{P} \Big \{ \delta \in \Delta :~ g( x_j,\delta) >0\Big \} \\
& \leq \sum_{j \in I_{d+1}}\max_{j \in I_{d+1}} \mathbb{P} \Big \{ \delta \in \Delta :~g( x_j,\delta) >0\Big \}  \\
&= (d+1) \max_{j \in I_{d+1}} \mathbb{P} \Big \{ \delta \in \Delta :~ g( x_j,\delta) >0\Big \} \\
&\, \leq (d+1) \max_{j=1,...,K_N} \mathbb{P} \Big \{ \delta \in \Delta :~g( x_j,\delta) >0\Big \} = (d+1)\ V( \hat{x}), 
\end{align*}
where  $V(\hat{x})=\max_{j=1,...,K_N} \mathbb{P} \Big \{ \delta \in \Delta :~g( x_j,\delta) >0\Big \}$ is the maximum probability of violation among the vertices $x \in v(\Pi_N)$. This concludes our proof. \qed \\
The inequality $V(x) \leq (d+1)\ V(\hat{x}), \ \forall \ x \in \Pi_N $ can be equivalently stated as:
There exists $\hat{x} \in v(\Pi_N)$ such that $\mathbb{V}(\Pi_N)=\sup_{x \in\Pi_N} V(x) \leq (d+1)V(\hat{x})$, which means that to bound the  probability of violation of the entire polytope, as defined in Definition \ref{def2}, we only need to know the number of decision variables and the vertex with the highest probability of violation. This relation forms a bridge between the notion of violation of a point and that of a set. 
\bibliography{biblio}             
 \setcitestyle{authoryear,open={((},close={))}}                                         
                                                   
\end{document}